\documentclass[12pt]{amsart}
\usepackage{amsmath}
\usepackage{amsrefs}\usepackage{amssymb}

\usepackage{dingbat}
\usepackage{skak}
\usepackage{amsfonts}\usepackage{verbatim}
\usepackage{marvosym}\usepackage{wasysym}
\usepackage{dictsym}
\usepackage{linearb}
\usepackage{amsbsy}
\usepackage{amsopn}\usepackage{MnSymbol}

\def\proclaim#1{\vskip0.5em\noindent{\bf #1}\it}
\def\endproclaim{\vskip0.5em\par\noindent\rm}
\def\proclaim#1{\vskip0.5em\noindent{\bf #1}\it}
\def\endproclaim{\vskip0.5em\par\noindent\rm}
\def\demo#1{\vskip0.5em\noindent{\bf #1\ }}

\def\text#1{\mbox{#1}}
\def\flushpar{\par\noindent}

\def\tag#1{\eqno{(#1)}}
\def\mod{\mbox{ mod }}
\newcommand{\mapright}[1]{%
    \smash{\mathop{%
        \hbox to 1cm{\rightarrowfill}
        }
    \limits^{#1}
    }
}
\newcommand{\mapleft}[1]{%
    \smash{\mathop{%
        \hbox to 1cm{\rightarrowfill}
        }
    \limits_{#1}
    }
}

\def\e{\epsilon}
\def\a{\alpha}
\def\b{\beta}

\def\D{\Delta}
\def\s{\sigma}
\def\Si{\Sigma}

\def\l{\lambda}
\def\x{\times}

\def\o{\overline}
\def\f{\flushpar}
\def\u{\underline}
\def\v{\varphi}

\def\om{\omega}
\def\Om{\Omega}
\def\B{\mathcal B}
\def\T{\widehat T}
\def\({\left(}
\def\){\right)}

\def\<{\langle}
\def\>{\rangle}
\def\bul{\smallskip\f$\bullet\ \ \ $}

\def\lfl{\lfloor}\def\rfl{\rfloor}\def\sms{\smallskip\f}
\def\pf{\smallskip\f{\it Proof}\ \ \ \ }\def\sms{\smallskip\f}\def\sbul{\f$\bullet\
\ \ $}\def\sms{\smallskip\f}
\def\Lra{\Longrightarrow}
\def\lfl{\lfloor}\def\rfl{\rfloor}\def\xbm{(X,\B,m)}\def\xbmt{(X,\B,m,T)}\def\xyr{\xrightarrow}
\def\undersetbrace#1\to#2{\underbrace{#2}_{#1}}

\def\({\biggl(}
\def\){\biggr)}
\def\<{\langle}
\def\>{\rangle}

\def\bul{\smallskip\f$\bullet\ \ \ $}\def\lfl{\lfloor}\def\rfl{\rfloor}

\def\){\biggr)}

\def\<{\bold\langle}
\def\>{\bold\rangle}

\def\bul{\smallskip\f$\bullet\ \ \ $}\def\sms{\smallskip\f}
\begin{document}\title{Symmetric Birkhoff sums in infinite ergodic theory   }
\author{ Jon. Aaronson}\thanks{Aaronson’s research was partially supported by ISF grant No. 1599/13.}
\author{Zemer Kosloff}\thanks{Kosloff's research was partially supported  by ERC AdG 320977.}
\author{Benjamin Weiss}
\address[Aaronson]{School of Math. Sciences, Tel Aviv University,
69978 Tel Aviv, Israel.}
\email{aaro@tau.ac.il}
\address[Kosloff]{Mathematics institute, University of Warwick, 
Coventry, CV4 7AL, UK.}
\email{z.kosloff@warwick.ac.uk}
\address[Weiss]{Institute of Mathematics
    Hebrew Univ. of Jerusalem,
    Jerusalem 91904, Israel}
\email{weiss@math.huji.ac.il}
\begin{abstract} We show that the  absolutely normalized, symmetric Birkhoff sums of positive integrable functions in infinite, 
ergodic systems never converge pointwise even though they may be almost surely bounded away from zero and infinity. We also consider the latter phenomenon characterizing it among
transformations admitting generalized recurrent events \end{abstract}\thanks{\copyright 2013++}
\maketitle
\section*{\S0 Introduction}
\subsection*{Pointwise ergodicity}
Pointwise ergodicity for infinite measure preserving transformations fails.
Let $\xbmt$ be a   conservative, ergodic , measure preserving transformation with $m(X)=\infty$ then  (see
\cite{IET}) for any sequence of constants $a_n>0$

\begin{align*}\tag*{\Large{$\diamondsuit$}}&\text{either}\ \varliminf_{n\to\infty}\frac{S_n(1_A)}{a_n}=0\ \ \text{a.e.}\ \ \forall\ A\in\B,\ m(A)<\infty\\ & \text{or}\ \exists\ n_k\to\infty\ \text{so that}\ \  \frac{S_{n_k}(1_A)}{a_{n_k}}\xyr[k\to\infty]{}\infty\ \ \text{a.e.}\ \ \forall\ A\in\B,\ m(A)>0.\end{align*}
This is for the one-sided Birkhoff sums $S_n(1_A)(x):=\sum_{k=0}^{n-1}1_A(T^kx)$.
\

For an invertible  $\xbmt$, and the one-sided Birkhoff sums replaced by the two-sided Birkhoff sums
$$\Si_n(1_A)(x):=\sum_{|k|< n}1_A(T^kx),$$ the $\diamondsuit$\  analogue may fail.
\

Infinite measure examples were given in \cite{MS} with constants $a_n>0$ so that
\begin{align*}\tag*{{\Large\Pointinghand}}\Si_n(1_A)\ \asymp\ a_n\ \text{a.e.} \ \forall\ \ A\in\B,\ 0<m(A)<\infty\end{align*}
where for eventually positive sequences $a_n,\ b_n$,   $a_n\ \asymp\ b_n$ means existence of $M>1$ so that
$M^{-1}<\tfrac{a_n}{b_n}<M\ \forall\ n$ large.

\

An early hint of this possibility can be found in \cite{DE} where an example $\xbmt$ is given for which the forward sums are not comparable to the backward sums, namely: for $A\in\B,\ 0<m(A)<\infty$ and a.e. $x\in X$,
\begin{align}
\varliminf_{n\to\infty}\frac{S^-_n(1_A)(x)}{S^+_n(1_A)(x)}=0\ \ \&\ \ \varlimsup_{n\to\infty}\frac{S^-_n(1_A)(x)}{S^+_n(1_A)(x)}=\infty \label{eq: rightthumbsdown}
\end{align}
where $S^+ _n(f)=S_n(f)$ and $S^-_n(f):=\sum_{k=0}^{n-1}f\circ T^{-k}$. Indeed,  in view of $\diamondsuit$,\ \  \eqref{eq: rightthumbsdown} \ \  is a consequence of {\Large\Pointinghand}, which is in turn satisfied by the \cite{DE} example (see theorem 3 below).

\

Our main result (theorem 2 in \S2) is that {\Large\Pointinghand} can never be  upgraded to the convergence:
\begin{align*}\tag*{{\leftthumbsup}}\frac{\Si_n(1_A)}{a_n}\ \xyr[n\to\infty]{}m(A)\ \ \text{a.e.}\ \ \ \forall\ \ A\in\B,\ 0<m(A)<\infty.\end{align*}

\
For an infinite measure preserving transformation, {\Large\Pointinghand} implies {\tt bounded rational ergodicity} (see \S1). We proceed to discuss two classes of bounded rationally ergodic transformations. In section 3 we show that rank one transformations with a bounded cutting sequence always satisfy {\Large\Pointinghand}. 
\

In section 4 we first prove Theorem 4 which states that the return sequence of a bounded rationally ergodic, weakly pointwise dual ergodic transformation satisfies the {\tt extended regular variation} property with Karamata indices 1, meanning that there exists $M\geq1$ and $N:\mathbb{N}\to\mathbb{N}$ so that
$$M^{-1}pa(n)\leq a(pn)\leq M pa(n),\ \forall p>1,\ n\geq N(p).$$
This is used to prove  a necessary and sufficient condition (Theorem 5 in \S4) for a transformation admitting a generalized recurrent event to satisfy {\Large\Pointinghand}. This condition is a "trimmed sum" type small tail condition of the first return time functions to generalized recurrent events.

For larger groups the situation is different. Examples of actions of large groups satisfying {\leftthumbsup} analogues for compactly supported, integrable  functions are given in  Theorem 1.1 of \cite{GN}. In this context (example 5.1 in \S5) we show that certain infinite $\Bbb Z^2$ actions satisfy the  {\leftthumbsup} analogue for all integrable functions. We also discuss the various possibilities for infinite ergodic $\Bbb Z^2$ actions in terms of the recurrence of the generators.

\
The ratio ergodic theorem holds for $\Bbb Z^d$ actions (see \cite{EH} for $d=1$ and \cite{H} for $d\ge 2$) and shows that for an ergodic $\Bbb Z^d$ action, 
if any of these statements holds for some $A\in\B,\ 0<m(A)<\infty$, then it holds for all $ A\in\B,\ 0<m(A)<\infty$. 
Thus  the properties are invariant under similarity (see \cite{IET}).
\

\section*{\S1 Preliminaries}
\subsection*{Bounded rational ergodicity}
\

As in \cite{BRE}, the  conservative, ergodic, measure preserving transformation  $\xbmt$ is called
{\it boundedly rationally ergodic} ({\tt BRE}) if $\exists\ A\in\B,\ 0<m(A)<\infty$ and $M>0$ so that
\begin{equation}
\label{eq: a}
\begin{aligned}&S_n(1_A)(x)\le Ma_n(A)\ \text{a.e. on}\ A\ \forall\ n\ge 1\\ & \text{where}\ a_n(A):=\sum_{k=0}^{n-1}\frac{m(A\cap T^{-k}A)}{m(A)^2}.\end{aligned}
\end{equation}
In this case (\cite{BRE}), $\xbmt$ is {\it weakly rationally ergodic} ({\tt WRE}), that is, writing $a_n(T):=a_n(A)$ (where $A$ is as in \eqref{eq: a}), there is a dense hereditary ring
$$R(T)\subset\mathcal F:=\{F\in\B:\ m(F)<\infty\}$$ (including all sets satisfying \eqref{eq: a}) so that 
$$\sum_{k=0}^{n-1}m(F\cap T^{-k}G)\sim m(F)m(G)a_n(T)\ \forall\ F,\ G\in R(T)$$
and in particular,
$$a_n(F)\sim a_n(T)\ \forall\ F\in R(T),\ m(F)>0.$$
\

For invertible transformations, the  one sided properties ({\tt RE} $\&$ {\tt BRE}) are equivalent to their 2-sided analogues:
$\xbmt$ is:
\bul {\it two-sidedly, boundedly rationally ergodic} if $\exists\ A\in\B,\ 0<m(A)<\infty$  and $M>0$ so that
\begin{equation} \label{eq: a'}
\begin{aligned}&\Si_n(1_A)(x)\le M\o a_n(A)\ \text{a.e. on}\ A\ \forall\ n\ge 1\\ & \text{where}\ \o a_n(A):=\sum_{k=-(n-1)}^{n-1}\frac{m(A\cap T^{k}A)}{m(A)^2}\sim 2a_n(A);\end{aligned}
\end{equation}
and
\bul {\it two-sidedly,  weakly rationally ergodic}, if    there is a dense hereditary ring
$$\o R(T)\subset\mathcal F$$ (including all sets satisfying \eqref{eq: a'}) so that $$\o a_n(F)\sim 2a_n(T)\ \forall\ F\in \o R(T),\ m(F)>0.$$

In case $T$ is weakly rationally ergodic, $\exists\ \u\b(T)\in [0,1],\ \a(T),\ \b(T)\in [1,\infty]$ so that a.e., $\forall\ f\in\ L^1(m)_+$:
\begin{align*}&\varlimsup_{n\to\infty}\frac1{a_n(T)}S_n(f)\ =\ \a\int_Xfdm\ \\ &
\varlimsup_{n\to\infty}\frac1{2a_n(T)}\Si_n(f)\ =\ \b\int_Xfdm\\ &
\varliminf_{n\to\infty}\frac1{2a_n(T)}\Si_n(f)\ =\ \u\b\int_Xfdm
\end{align*}
and $T$ is boundedly rationally ergodic if and only if $\a(T)<\infty$. See \cite{BRE} for the one sided Birkhoff sums, the case of symmetric Birkhoff sums is similar. 

Since bounded rational ergodicity of an invertible $\xbmt$ implies that of $T^{-1}$ we have that $\a(T)<\infty\ \Lra\ \a(T^{-1})<\infty$. Thus an invertible, bounded rationally ergodic transformation is two-sidedly bounded rationally ergodic. 

\proclaim{Proposition 1} \label{prop: connection between b and a}
\

\  Let $\xbmt$ be an invertible, conservative, ergodic , measure preserving transformation.
  \sms{\rm (i)} If \ \ $T$ satisfies {\Large\Pointinghand} wrt some sequence of normalizing constants, then $T$ is boundedly rationally ergodic, (hence
  weakly rationally ergodic).

\sms{\rm (ii)} If \ \ $T$ is  boundedly rationally ergodic, then

\begin{align}&\label{eq: clubsuit}\b(T)\le \a(T)=\a(T^{-1})\le 2\b(T)\ \&\\ & \label{eq: Stopsign}
\ \ \u\b(T)\ \le\ \frac{\a(T)}2.\end{align}\endproclaim\demo{Proof of {\rm (i) }}\ \ \ \ \ Suppose that 
$$\Si_n(1_A)\ \asymp\ \ a(n)m(A)\ \text{a.e}\ \text{for some and hence all}\ A\in\B,\ 0<m(A)<\infty.$$
Fix $A\in\B,\ 0<m(A)<\infty$. By Egorov's theorem $\exists\ M>1,\ N\in\Bbb N\ \&\ B\in\B(A),\ m(B)>0$ so that
$$\frac{\Si_n(1_A)(x)}{a(n)}<Mm(A)\ \ \forall\ x\in B,\ n\ge N.$$
On the other hand, $\exists\ \e>0$ so that
$$\varliminf_{n\to\infty}\frac{\Si_n(1_B)}{a(n)}\ge 4\e \ \text{a.e.}$$
whence, by Fatou's lemma
$$a_n(B)\ge\frac{1}{3m(B)^2}\int_B\Si_n(1_B)dm\ge\e \frac{a(n)}{m(B)}\ \forall\ n\ \ \text{large}.$$
To see bounded rational ergodicity, for $n\ge 1$ large and $x\in B$,
$$S_n(1_B)(x)\le \Si_n(1_A)(x)\le Ma(n)m(A)\le\frac{Mm(A)m(B)}\e\cdot a_n(B).\ \ \ \CheckedBox\text{(i)}$$
\demo{Proof of {\rm(ii) \eqref{eq: clubsuit}}}
\

It suffices to show that $\a(T)\ge\a(T^{-1})$. Fix $A\in\mathcal F_+\ \&\ \e>0$. By Egorov, $\exists\ B\in\B_+\cap A\ \&\ N_0\ge 1$ so that
$$S_n^+(1_B)(x)<(\a+\e)a_n(T)m(B)\ \forall\ x\in B,\ n\ge N_0.$$
For $n\ge 1\ \&\ x\in B$, let
$$K_n(x):=\max\,\{0\le k\le n:\ T^{-k}(x)\in B\};$$
then $K_n\xyr[n\to\infty]{}\infty$ a.e. $\&$ whenever $K_n(x)\ge N_0$,
\begin{align*}
    S_n^-(1_B)(x)&=S^-_{K_n(x)}(1_B)=S_{K_n(x)}(1_B)\circ T^{-K_n}(x)\\ &\leq (\a+\e)a_{K_n(x)}(T)m(B)\\ & \le (\a+\e)a_n(T)m(B).\ \CheckedBox\ \eqref{eq: clubsuit}
\end{align*}

\demo{Proof of {\rm (ii) \eqref{eq: Stopsign}}}\
\
Let $B\in\mathcal{F}$. Since $T$ is bounded rationally ergodic, $\alpha(T)=\alpha(T^{-1})<\infty$ and thus 
by $\diamondsuit$, 
$$\varliminf_{n\to\infty} \frac{S^-_n(1_B)}{a_n(T)}=0\ \ \text{a.e}.$$
It follows that a.e.,
\begin{align*}
\varliminf_{n\to\infty} \frac{\Si_n(1_B)}{2a_n(T)}&=\varliminf_{n\to\infty} \frac{S^+_n(1_B)+S^-_n(1_B)}{2a_n(T)}\\ &\leq \varliminf_{n\to\infty} \frac{S^-_n(1_B)}{2a_n(T)}+\varlimsup_{n\to\infty} \frac{S^+_n(1_B)}{2a_n(T)}\\
 &=\frac{\alpha}{2} m(B).
\end{align*}


\section*{\S2 No absolutely normalized convergence of two-sided Birkhoff sums.}

 \proclaim{Theorem 2}
 \

  \  Let $\xbmt$ be an infinite, invertible, conservative, ergodic , measure preserving transformation, then
 {\leftthumbsup} fails.
 \endproclaim\demo{Proof}\  \ Suppose otherwise, namely that for some $a(n)>0$,
\begin{align*}\tag*{{\leftthumbsup}}\frac{\Si_n(1_A)}{2a(n)}\ \xyr[n\to\infty]{}m(A)\ \ \text{a.e.}\ \ \forall\ \ A\in\B,\ 0<m(A)<\infty\end{align*}
By proposition 1(i), $T$ is boundedly rationally ergodic, hence weakly rationally ergodic. 
\

We claim first that  $a(n)\sim a_n(T)$. To see this, let $A\in R(T)$ and let $B\in\B(A),\ m(B)>0$ so that
$$\frac{\Si_n(1_A)}{2a(n)}\ \xyr[n\to\infty]{}m(A)\ \ \text{uniformly on}\ B.$$
It follows that
$$\sum_{k=0}^{n-1}m(B\cap T^{-k}A)\sim\frac12\int_B\Si_n(1_A)dm\sim m(A)m(B)a(n).$$
On the other hand, since $A\in R(T)$ and $B\subset A$, by Theorem 3.3.1 in \cite{IET},
$$\sum_{k=0}^{n-1}m(B\cap T^{-k}A)\sim m(A)m(B)a_n(T)$$
showing that indeed  $a(n)\sim a_n(T)$.

\

We claim next that $\a(T)=2$.

\

Indeed by \eqref{eq: Stopsign}, $\a\ge 2$ and by
{\eqref{eq: clubsuit}}, $\a\le 2$. Thus

$$\varlimsup_{n\to\infty}\frac1{2a(n)}S_n(1_F)=m(F)\ \text{a.e.}\ \ \forall\ F\in\B.$$
\

The rest of  this proof is on a ``single orbit" which  we proceed to specify.
\

\bul Fix  $A\in\mathcal F_+$. By Egorov, $\exists\ B\in\B(A),\ m(B)>\frac34m(A)$ so that
$$\sup_{N\ge n}\frac1{2a(N)}S_N(1_A),\ \frac1{2a(n)}\Si_n(1_A)\xyr[n\to\infty]{}\ m(A)\ \text{uniformly on}\ B.$$

\bul Call a point $x\in B$ {\it admissible} if
\begin{align*}&\tag*{A(i)}\frac{S_N(1_B)(x)}{S_N(1_A)(x)}\ \xyr[n\to\infty]{}\
\frac{m(B)}{m(A)};\\ &\tag*{A(ii)} \frac1{2a(n)}\Si_n(1_B)(x)\xyr[n\to\infty]{}\ m(B)\\ &\tag*{A(iii)}\sup_{N\ge n}\frac1{2a(N)}S_N(1_B)\xyr[n\to\infty]{}\ m(B),\
\end{align*}
and $\exists$  $K\subset\Bbb N$, an {\it $x$-admissible subsequence} in the sense that
\begin{align*}&\tag*{A(iv)} T^nx\in B\ \forall\ n\in K\ \ \&\ \\ &\tag*{A(v)} \frac1{2a(n)}S_n(1_B)(x)\xyr[n\to\infty,\ n\in K]{}\ m(B).
\end{align*}
An {\it admissible pair} is $(x,K)\in B\x 2^\Bbb N$ where $x$ is an {admissible point} and $K$ is an {$x$-admissible subsequence}.
\

Note that if $(x,K)$ is an admissible pair, then by A(iv) and A(i),
$$\frac1{2a(n)}S_n(1_A)(x)\xyr[n\to\infty,\ n\in K]{}\ m(A).$$
In what follows $c_n\lesssim d_n$ means $\varlimsup_{n\to\infty}\frac{c_n}{d_n}\leq 1$.
\proclaim{Lemma 0}\ \ Almost every $x\in B$ is admissible.\endproclaim\demo{Proof}
\

By {\leftthumbsup}, $\a(T)=2$ and the ratio theorem, almost every $x\in B$ satisfies A(i), A(ii) $\&$ A(iii).
\

Also by  $\a(T)=2$, for a.e. $x\in B,\ \exists\ K\subset\Bbb N$ satisfying A(v).

\

We claim that if $K:=\{k_n:\ n\ge 1\},\ k_n\uparrow$, then
$K':=\{k'_n:\ n\ge 1\}$ where $k'_n:=\max\{j\le k_n:\ T^jx\in B\}$ is $x$-admissible. Evidently $K'$  is infinite and satisfies  A(iv). To check A(v):
$$2a(k_n)m(B)\ge 2a(k_n')m(B)\overset{\text{\tiny A(iii)}}{\text{\Large$\gtrsim$}} S_{k_n'}(1_B)(x)=S_{k_n}(1_B)(x)\overset{\text{\tiny A(v)}}{\text{\Large$\sim$}}2a(k_n)m(B).\ \ \CheckedBox\text{A(v)}$$
\proclaim{Lemma 1}\ \ \ If $x\in B,\ K\subset\Bbb N$ and $\{J_n:\ n\in K\}$ satisfy
\begin{align*}&\frac1{2a(n)}S_n(1_A)(x)\xyr[n\to\infty,\ n\in K]{}\ m(A);
\\ &\ n\ge J_n\xyr[n\to\infty,\ n\in K]{}\ \infty;\\ &\ \varliminf_{n\to\infty,\ n\in K}\frac{a(J_n)}{a(n)}=:\rho>0,
\end{align*}
then
$$\frac1{2a(J_n)}S_{J_n}(1_A)(x)\xyr[n\to\infty,\ n\in K]{}\ m(A).$$\endproclaim

\demo{Proof}
\begin{align*}
 \frac1{2a(J_n)}S^{-}_{J_n}(1_A)(x) &\lesssim \frac1{2\rho a(n)}S^{-}_{n}(1_A)(x)\ \ \text{as}\ n\to\infty,\ n\in K;\\ &=
 \frac1{\rho}\left(\frac1{2a(n)}\Si_{n}(1_A)(x)-\frac1{2 a(n)}S_{n}(1_A)(x)\right)\\ &\xyr[n\to\infty,\ n\in K]{}\ 0\ \ \ \ \because\ \ x\in B.
\end{align*}
Therefore,
\begin{align*}\ \ \frac1{2a(J_n)}S_{J_n}(1_A)(x) &= \frac1{2a(J_n)}\Si_{J_n}(1_A)(x)- \frac1{2a(J_n)}S^{-}_{J_n}(1_A)(x)\\ &\xyr[n\to\infty,\ n\in K]{}\ m(A).\ \ \CheckedBox\end{align*}
\proclaim{Lemma 2}\ \ \ Let $(x,K)\in B\x 2^\Bbb N$ be an admissible pair, then
\begin{align*}\frac1{12}\le \varliminf_{n\to\infty,\ n\in K}\frac{a(\frac{n}9)}{a(n)}\ \ \&\ \  \ \ \varlimsup_{n\to\infty,\ n\in K}\frac{a(\frac{n}9)}{a(n)}\le\frac14.
\end{align*}
\endproclaim
\demo{Proof}\ \ We show first that
\begin{align}\label{eq: a1}\varliminf_{n\to\infty,\ n\in K}\frac{a(\frac{n}9)}{a(n)}\ge\frac1{12}.
\end{align}

\ \ \ Define
$$J_i:=\min\,\{\ell\ge \frac{in}9:\ T^\ell x\in B\}\wedge\frac{n(i+1)}9;\ \ \ (0\le i\le 8),$$
then
\begin{align*}2m(B)a(n)&\lesssim\ S_n(1_B)(x)\ \ \text{as}\ n\to\infty,\ n\in K\\ &=\sum_{i=0}^8S_{\frac{n}9}(1_B)(T^{\frac{in}9}x)\\ &=
\sum_{i=0}^8S_{\frac{(i+1)n}9-J_i}(1_B)(T^{J_i}x)
\\ &\le \sum_{i=0}^8S_{\frac{n}9}(1_B)(T^{J_i}x)\\ &\le \sum_{i=0}^8\|S_{\frac{n}9}(1_B)\|_{L^\infty(B)}
\ \ \because\ \|S_{\frac{n}9}(1_B)\|_{L^\infty(X)}=\|S_{\frac{n}9}(1_B)\|_{L^\infty(B)}
\\ &\lesssim 18m(A)a(\frac{n}9)\ \ \text{as}\ n\to\infty.
\end{align*}
Thus
$$\varliminf_{n\to\infty,\ n\in K}\frac{a(\frac{n}9)}{a(n)}\ge\frac{2m(B)}{18m(A)}>\frac1{12}.\ \ \ \CheckedBox{\eqref{eq: a1}}$$
Next, we show:
\begin{align}\label{eq: b}\varlimsup_{n\to\infty,\ n\in K}\frac{a(\frac{n}3)}{a(n)}\le\frac12.
\end{align}

By  \eqref{eq: a1}, $\{\frac{n}3:\ n\in K\}$ satisfies the preconditions of Lemma 1 and so
\begin{align}\label{eq: Mortiss} \frac1{2a(\tfrac{n}3)}S_{\frac{n}3}(1_A)(x)\ \xyr[n\to\infty,\ n\in K]{}\ m(A).
\end{align}

For $n\in K$, let
$$J_n:=\max\,\{j< \frac{n}3:\ T^jx\in B\}.$$
We claim that
$a(J_n)\sim a(\tfrac{n}3)$ as $n\to\infty,\ n\in K$ since:
\begin{align*}2a(\tfrac{n}3)m(B)&\ge 2a(J_n)m(B)\\ &\underset{n\to\infty}{\text{\Large$\gtrsim$}} S_{J_n}(1_B)(x)\ \ \because\ x\in B;\\ &=S_{\frac{n}3}(1_B)(x)\\ &\underset{n\to\infty,\ n\in K}{\text{\Large$\sim$}}2a(\tfrac{n}3)m(B).
\end{align*}
The last step uses A(i) and \eqref{eq: Mortiss}. 

Finally as $n\to\infty,\ n\in K$
\begin{align*}2m(A)a(n)&\sim\ \Si_n(1_{A})(T^{J_n}x)\ \ \because\ T^{J_n}x\in B;\\ &=\sum_{k=-n+J_n}^{n+J_n}1_A(T^kx)\\ &\ge
\Si_{J_n}(1_{A})(T^{J_n}x)+\Si_{J_n}(1_A)(T^nx)\\ &\sim 4m(A)a(J_n)\ \ \ \ \because\ T^{J_n}x,\ T^nx\in B;\\ &\sim 4m(A)a(\tfrac{n}3).\ \ \CheckedBox{\rm \eqref{eq: b}}
\end{align*}
Next, we iterate \eqref{eq: b}:
\begin{align}\label{eq: c}\varlimsup_{n\to\infty,\ n\in K}\frac{a(\frac{n}9)}{a(n)}\le\frac14.
\end{align}

\demo{Proof of \eqref{eq: c}}\ \
Let $L_n:=\min\,\{J\ge\frac{n}3:\ T^Jx\in B\}$.
\

We claim that
\begin{align}\label{eq: Lightning} \exists\ N\ge 1\ \ \text{so that}\  L:=\{L_n:\ n\in K\ n\ge N\}\ \ \text{ is $x$-admissible}.
\end{align}

\demo{Proof}
\

Firstly, for $n\in K,\ T^nx\in B$ whence  $L_n\le n$.

Since
$$\varliminf_{n\to\infty,\ n\in K}\frac{a(L_n)}{a(n)}\ge \varliminf_{n\to\infty,\ n\in K}\frac{a(\frac{n}3)}{a(n)}\ \ \overset{\eqref{eq: a1}}{\text{\Large$\ge$}}\ \ \frac1{12}$$
we have by Lemma 1 and A(i) that
$$\frac1{2a(L_n)}S_{L_n}(1_B)(x)\xyr[n\to\infty,\ n\in K]{}\ m(B)$$
and
$L:=\{L_n:\ n\in K\}$ is $x$-admissible.\ \Checkedbox\text{\eqref{eq: Lightning}}

\

Applying \eqref{eq: b}  to $L=\{L_n:\ n\in K\}$, we obtain
\begin{align*}\varlimsup_{n\to\infty,\ n\in K}\frac{a(\frac{L_n}3)}{a(L_n)}\le \frac12.\end{align*}
Since $a(n/9)\leq a(L_n/3)$, to obtain \eqref{eq: c} from this, it suffices to show that
\begin{align}\label{eq: Biohazard} a(L_n)\underset{n\to\infty,\ n\in K}{\text{\Large $\sim$}}\ a(\frac{n}3)
\end{align}
\demo{Proof of\ \eqref{eq: Biohazard}}
\

By A(i) and \eqref{eq: Mortiss},
$$\frac1{2a(\frac{n}3)}S_{\frac{n}3}(1_B)(x)\xyr[n\to\infty,\ n\in K]{}\ m(B)$$
whence, as $n\to\infty,\ n\in K$
\begin{align*}2a(\frac{n}3)m(B)&\ \sim\ S_{\frac{n}3}(1_B)(x)\\ &\sim S_{L_n}(1_B)(x)\\ &\sim 2a(L_n)m(B)\\ &\ge 2a(\frac{n}3)m(B).\ \CheckedBox\text{\
\eqref{eq: Biohazard}\ \&\ \eqref{eq: c}}
\end{align*}
This  completes the proof of lemma 2.\ \ \ \Checkedbox

\begin{align*}\tag*{\bf Lemma 3:}\ \text{If}\ & (x,K) \ \text{is an admissable pair},\ \l,\ \rho\in (0,1)\ \&\ \varliminf_{n\to\infty,\ n\in K}\frac{a(\l n)}{a(n)}\ge\rho,\ \ \text{then}\\ &
\varlimsup_{n\to\infty,\ n\in K}\frac{a((1-\l) n)}{a(n)}\le\ 1-\rho.
 \end{align*}
 \demo{Proof}
 \

 Firstly we claim that as $n\to\infty,\ n\in K$,
 \begin{align}\label{eq: Coffeecup}S^{-}_{\l n}(1_A)(T^nx)\sim 2m(A)a(\l n).
 \end{align}

 To see \eqref{eq: Coffeecup}, note that as $n\to\infty,\ n\in K$,
 $$S^{-}_n(1_A)(T^nx)=S_n(1_A)(x)\sim 2m(A)a(n).$$
 Since $T^nx\in B$, we have that
 \begin{align*}\frac1{2a(n)}S_n(1_A)(T^nx)&=\frac1{2a(n)}\Si_n(1_A)(T^nx)-
 \frac1{2a(n)}S^{-}_n(1_A)(T^nx)\\ & \xyr[n\to\infty,\ n\in K]{}\ 0
 \end{align*}
 whence also
 $$\frac1{2a(\l n)}S_{\l n}(1_A)(T^nx)\ \xyr[n\to\infty,\ n\in K]{}\ 0$$ and
 \begin{align*}\frac1{2a(\l n)}S^{-}_{\l n}(1_A)(T^nx)&=\frac1{2a(\l n)}\Si_{\l n}(1_A)(T^nx)-\frac1{2a(\l n)}S^{+}_{\l n}(1_A)(T^nx)\\ & \xyr[n\to\infty,\ n\in K]{}\ m(A).\ \CheckedBox{\eqref{eq: Coffeecup}}
 \end{align*}

 To prove the lemma, we  assume without loss of generality,  that there exists $\varepsilon>0$ such that 
 $$\varlimsup_{n\to\infty,\ n\in K}\frac{a((1-\l) n)}{a(n)}>\varepsilon$$
 (if this fails the lemma holds) and thus by Lemma 1,
 $$\frac{1}{2a((1-\l)n)}S_{(1-\l)n}\left(1_A \right) \xrightarrow[n\to\infty,\ n\in K]{}m(A).$$
 Consequently as $n\to\infty,\ n\in K$,
 \begin{align*}2m(A)a((1-\l)n)&\sim S_{(1-\l)n}(1_A)(x)\\ &=
 S_{n}(1_A)(x)-S^{-}_{\l n}(1_A)(T^nx)\\ &\sim 2m(A)(a(n)-a(\l n))\\ &\lesssim 2(1-\rho)m(A)a(n). \ \CheckedBox
 \end{align*}
 \demo{Proof of the theorem}
 \

 Fix an admissible pair $(x,K)\in B\x 2^\Bbb N$, then
by Lemmas 2 and 3 $$\varlimsup_{n\to\infty,\ n\in K}\frac{a(\frac{n}9)}{a(n)}\le\frac14\ \ \&\ \ \ \varlimsup_{n\to\infty,\ n\in K}\frac{a(\frac{8n}9)}{a(n)}\le\frac{11}{12}.$$
 For $n\in K$, let
 $$J_n=J_n(x):=\min\,\{j\ge \frac{n}9:\ T^jx\in B\}.$$
 We claim that $J_n\le\frac{8n}{9}$; else, as $n\to\infty,\ n\in K$:
 \begin{align*}2a(n)m(B)&\lesssim S_n(1_B)(x)\\ &=
 S_{J_n}(1_B)(x)+S_{(n-J_n)\vee 0}(1_B)(T^{J_n}x)\\ &=
  S_{\frac{n}9}(1_B)(x)+S_{(n-J)\vee 0}(1_B)(T^{J_n}x)
    \\ &\le S_{\frac{n}9}(1_B)(x)+S_{\frac{n}9}(1_B)(T^{J_n}x)\ \text{assuming}\ J_n>\frac{8n}{9};\\ &\lesssim 4m(A)a(\frac{n}9)
 \end{align*}
 whence as $n\to\infty$, $n\in K$,
 $$\frac{a(\frac{n}9)}{a(n)}\gtrsim \frac{2m(B)}{4m(A)}>\frac{3}{8}$$
 and
 $$\frac38\le \varlimsup_{n\to\infty,\ n\in K}\frac{a(\frac{n}9)}{a(n)}\le\frac14.\ \ \ \ \XBox$$
 This contradiction shows that indeed for sufficiently large $n\in K$, $J_n\le\frac{8n}{9}$.

 \

 Finally,  since $\frac{n}9\le J_n\le \frac{8n}9$: $[J_n-\frac{8n}9,J_n+\frac{8n}9]\supset [0,n]$ and as $n\to\infty,\ n\in K$,
 \begin{align*}
  2a\left(\frac{8n}9\right)m(A)&\sim \Si_{\frac{8n}9}(1_A)(T^{J_n}x)\\ &\ge S_n(1_A)(x)\\ &\sim 2a(n)m(A)
 \end{align*}
whence
$$\frac{11}{12}\ \underset{n\to\infty,\ n\in K}{\text{\Large$\gtrsim$}}\ \frac{a(\frac{8n}9)}{a(n)}\ \underset{n\to\infty,\ n\in K}{\text{\Large$\gtrsim$}}\ \ 1.\ \XBox$$
 This last contradiction shows that there is no ergodic theorem for symmetric Birkhoff sums of conservative ergodic infinite measure preserving transformations.\ \Checkedbox
 \subsection*{Remark on quantitative estimates}
 \

 The proof of theorem 2 can be adapted to show that $\exists\ \D>0$ so that for any $\xbmt$ satisfying {\Large\Pointinghand}, we have
 \begin{align*}\b(T)-\underline{\b}(T)\ \ge\ \D.
 \end{align*}
The question of estimating the best $\D>0$ arises. For the examples appearing in this paper, $\D\ge\frac12$.
 \section*{ \S3 Rank one towers}

These are {\tt CEMPT}s constructed by cutting and stacking as in
\cite{FR-Book}, \cite{FR-AMM},\ \cite{CH}, Ch. 7 of \cite{N2}.
\

Let $c_n\in \Bbb N,\ c_n\ge 2\ \ (n\ge 1)$ and let $S_{n,k}\ge 0,\ \ (n\ge 1,\ 1\le k\le c_n)$.
The {\it rank one transformation} with {\it construction data}
$$\{(c_n;S_{n,1},\dots,S_{n,c_n}):\ n\ge 1\}$$
is an invertible piecewise translation of the interval $J_T=(0,S_T)$ where
$$S_T:=1+\sum_{n\ge 1}\frac1{C_n}\sum_{k=1}^{c_n}S_{n,k}\le \infty\ \text{with}\ C_n:=c_1\cdots c_n.$$
This is defined as the limit of a
 nested sequence of Rokhlin towers $(\tau_n)_{n\ge 1}$ of intervals where $\tau_1=[0,1]$ and $\tau_{n+1}$ is constructed  from $\tau_n$ by
  \sbul    cutting $\tau_n$ into $c_n$ columns,\sbul putting $S_{n,k}$ spacer intervals above the $k^{\text{th}}$ column ($1\le k\le c_n$);\sbul  and stacking.

    \

     The transformation $T$ constructed, being an invertible,
    piecewise translation of $J_T$, preserves Lebesgue measure. It is  conservative and  ergodic.
       \

     \proclaim{Theorem 3} \label{thm: rank one}\ \ \ \ Let $\xbmt$ be the conservative, ergodic measure preserving transformation with construction data $$\{(c_n;S_{n,1},\dots,S_{n,c_n}):\ n\ge 1\}.$$
     \

     If $\sup_{n\ge 1}c_n<\infty$, then $T$ satisfies {\Large\Pointinghand}.\endproclaim\demo{Proof}\ \ Let $c_n\le J\ \ (n\ge 1)$ and let $q_n$ be the height of $\tau_n$  ($n\ge 1$).
     \

     For $x\in I:=[0,1]$ and $n\ge 1$ we have
     $$C_n\le\Si_{q_n}(1_I)(x)\le 2C_{n}.$$
     Define $a(n)$ by
     $$a(n):=C_\nu\ \ \text{for}\ \ q_\nu\le n<q_{\nu+1},$$
     then, for $q_\nu\le n<q_{\nu+1}$
     $$a(n)=C_\nu\le\Si_{q_\nu}(1_I)(x)\le \Si_{n}(1_I)(x)\le\Si_{q_{\nu+1}}(1_I)(x)\le 2C_{\nu+1}\le 2Ja(n).$$
     Finally $\varlimsup_{n\to\infty}\frac{\Si_n(1_A)}{a(n)}\ \&\ \varliminf_{n\to\infty}\frac{\Si_n(1_A)}{a(n)}$ are $T$-invariant whence constant by ergodicity and we have {\Large\Pointinghand}.\ \ \Checkedbox
\subsection*{Remarks}\ \ 
\

1.  By Theorem 3 and Proposition 1, every rank one transformation with bounded cuts and an arbitrary spacer sequence is bounded rationally ergodic. 
This was established also in \cite{DGPS}. 

2.  There are examples of bounded rationally ergodic rank one transformations whose return sequence grows arbitrarily slowly. Namely, if
$L_n\to\infty$, then there is a  rank one transformation  $T$  with cutting sequence $c_n\equiv 2$ with $a_n(T)=o(L_n)$ as $n\to\infty$.
See theorem 1 in \cite{BRE}. 
 
\section*{\S4 Weakly pointwise dual ergodic transformations}
As in \cite{AZ}, the conservative ergodic measure preserving transformation\ \ $(X,\B,m,T)$
is called {\it weakly pointwise dual ergodic}\ if
 $\exists\  a(n)>0, n\ge 1,$ such that for each $f\in L^1_+(m)$,
\begin{align*}&
\frac1{a(n)}\sum_{k=0}^{n-1}\T^k f \xrightarrow[n\to\infty]{m}\ \ 
\int_Xfdm\ \ \&\\ &
\varlimsup_{n\to\infty}\frac1{a(n)}\sum_{k=0}^{n-1}\T^k f=\int_Xfdm\ \ \text{a.e.}\ 
\end{align*}
where $\hat{T}:L_1(m)\to L_1(m)$ is the transfer operator defined by 
$$\forall f\in L_1(m), g\in L_\infty(m),\ \int_X \hat{T}f\cdot gdm=\int_X f\cdot g\circ Tdm.$$
This property entails {\tt WRE} and the return sequence $a_n(T)\sim a(n)$.
\

Our next result shows that the return sequence of bounded rationally ergodic, weakly pointwise dual ergodic transformation must be large. This is in contrast with the rank one transformations considered in theorem 3 whose return sequences can grow arbitrarily slowly.

\

\proclaim{Theorem 4}
\

Let \ $\xbmt$ be weakly  pointwise dual ergodic with return sequence $a(n)=a_n(T)$. 
\

If  $\a(T)<\infty$ then $\exists\ M>1$ and $N:\Bbb N\to\Bbb N$ so that
\begin{align*}\tag{{\tt ER}}a(pn)=M^{\pm 1}pa(n)\ \forall\ p>1,\ n\ge N(p).
\end{align*}
\endproclaim
\subsection*{Remark}\ \ The property ({\tt ER}) is called 
{\tt extended regular variation} with Karamata indices $1$  in  \cite{BGT}.

\demo{Proof of ({\tt ER})}
\

Fix  $\Om\in\B,\ m(\Om)=1$ a limited set in the sense of \cite{AZ}, that is satisfying
$$\left\|\frac1{a(n)}\sum_{k=0}^{n-1}\T^k 1_\Om\right\|_{L^\infty(\Om)}\xrightarrow[n\to\infty]{}\ 1.$$
 \

 WLOG, $a(n)=\sum_{k=0}^nu_k$ where $u_n:=m(\Om\cap T^{-n}\Om)$.
 
By lemma 4.1 in \cite{AZ},
\begin{align}&\label{eq: a2}\frac{\widehat{U}_s}{u(s)}\xrightarrow[s\to 0+]{m}1\ \ \&\\ &\label{eq: b1} \ \frac1{u(s)}\|\widehat{U}_s\|_{L^\infty(\Om)}\xrightarrow[s\to 0+]{}1
\end{align}
where $\widehat{U}_s:=\sum_{n\ge 0}e^{-sn}\widehat{T}^n1_{\Om}$ and $u(s):=\sum_{n\ge 0}e^{-sn}u_n$.

\

For $s>0$, set
$$U_s:=\sum_{n\ge 0}e^{-sn}1_{\Om}\circ T^n,$$
then $$\int_{\Om} U_sdm=\sum_{n\ge 0}u_ne^{-sn}=:u(s).$$
We claim first that for $p\in\Bbb N$

\begin{align}\label{eq: dsliterary}
 \int_\Om U_s^pdm\ \sim\ \ p!\prod_{k=1}^pu(ks)\ \ \text{as}\ s\to 0+.
\end{align}
\demo{\tt Proof of \eqref{eq: dsliterary}}
\

Firstly, by convexity, 
\begin{align}
 \label{eq: dsaeronautical}\int_\Om U_s^pdm\ \ge\ \(\int_\Om U_sdm\)^p=u(s)^p\ \forall\ p\ge 1,\ s>0.
\end{align}

Next, for $p\in\Bbb N$ fixed,

\begin{equation} \label{eq: dsarchitectural}
\begin{aligned}
 U_s^p=&p!V(p,s)+E(p,s)U_s^{p-1}\ \text{where}\\ & \ V(p,s)=\sum_{0\le n_1\le\dots\le n_p}e^{-s(n_1+\dots+n_p)}\prod_{k=1}^p1_{\Om}\circ T^{n_k}
 \\&\text{and}\ \ |E(p,s)|\le M_p\ \forall\ s>0\ \text{where $M_p$ is constant.}
\end{aligned}
\end{equation}
Thus

Now
\begin{align*}
 V(p,s)&=\sum_{n=0}^\infty e^{-ns}1_{\Om}\circ T^n\sum_{n\le n_2\le\dots\le n_p}e^{-s(n_2+\dots+n_p)}\prod_{k=2}^p1_{\Om}\circ T^{n_k}\\ &\overset{n_k=n+\nu_k}{\text{\Large$=$}}\
\sum_{n=0}^\infty e^{-ns}1_{\Om}\circ T^n\sum_{0\le \nu_2\le\dots\le \nu_p}e^{-s((p-1)n+\nu_2+\dots+\nu_p)}\prod_{k=2}^p1_{\Om}\circ T^{n+\nu_k}\\ &
=\sum_{n=0}^\infty e^{-nps}1_{\Om}\circ T^n\sum_{0\le \nu_2\le\dots\le \nu_p}e^{-s(\nu_2+\dots+\nu_p)}\prod_{k=2}^p1_{\Om}\circ T^{n+\nu_k}\\ &=
\sum_{n=0}^\infty e^{-nps}1_{\Om}\circ T^nV(p-1,s)\circ T^n.\end{align*}
whence
\begin{align*}\int_{\Om} V(p,s)dm&=\int_{\Om} V(p-1,s)\widehat{U}_{ps}dm
 \\ &\overset{\eqref{eq: b1}}{\underset{s\to 0+}{\text{\Large $\lesssim$}}}\ u(ps)\int_{\Om} V(p-1,s)dm.
\end{align*}

Thus
\begin{align}\label{eq: dsagricultural}
\int_{\Om} V(p,s)dm\underset{s\to 0+}{\text{\Large $\lesssim$}}\ \prod_{k=1}^pu(ks).
\end{align}
So far, by \eqref{eq: dsaeronautical},\ \eqref{eq: dsarchitectural} and \eqref{eq: dsagricultural},\ we have
\begin{align}\label{eq: dsmathematical} u(s)^p\ \asymp\ \int_\Om U_s^pdm\ =\ p!\int_\Om V(p,s)dm+O(u(s)^{p-1}). 
\end{align}
Thus, to finish the proof of \eqref{eq: dsliterary}, it suffices to show that
\begin{align}\label{eq: dsheraldical}\int_{\Om} V(p,s)dm\underset{s\to 0+}{\text{\Large $\gtrsim$}}\ u(ps)\int_{\Om} V(p-1,s)dm.
\end{align}
To this end, using \eqref{eq: dsmathematical} and \eqref{eq: dsagricultural}, we see that
  $$\int_{\Om} U_s^{2p}dm=O((\int_{\Om} U_s^{p}dm)^2)\ \text{as}\ s\to 0+$$
because
$$\int_{\Om} U_s^{2p}dm\asymp u(s)^{2p}=(\int_\Om U_sdm)^p\cdot(\int_\Om U_sdm)^p\overset{\text{\tiny convexity}}{\text{\Large $\le$}}(\int_\Om U_s^pdm)^2.$$
We'll need to know that
\begin{align}\label{eq: dschemical}
\int_{A}U_s^pdm\underset{s\to 0+}{\text{\Large $\sim$}}\ \ m(A)\int_\Om U_s^pdm\ \forall\ A\in\B(\Om). 
\end{align}
\demo{\tt Proof of \eqref{eq: dschemical}}
\

Let
{$$\Phi_s:=\frac{U_s^p}{\int_\Om U_s^pdm},$$} then $\int_\Om\Phi_sdm=1\ \&\ \sup_{s>0}\int_\Om\Phi_s^2dm<\infty$.
Thus $\{\Phi_s:\ s>0\}$ is weakly sequentially precompact in $L^2(\Om)$ and for \eqref{eq: dschemical}, it  suffices to show that
$$\Phi_s\xrightarrow[s\to 0]{}\ 1\ \ \text{weakly in}\ \ L^2(\Om).$$
To see this note that
\begin{align}\label{eq: ddag}
 e^{-s\v(x)}U_s\circ T_\Om(x)=U_s(x)-1 \ \text{on}\ \Om
\end{align}

where $\v(x):=\min\,\{n\ge 1:\ T^nx\in\Om\}$ (aka the {\it first return time} function) and $T_\Om$ is the {\it induced transformation} on $\Om$   defined by $T_\Om x:=T^{\v(x)}x$. 
\

As is well known,  $(\Om,\B(\Om),m_\Om,T_\Om)$ is an ergodic probability preserving transformation
 where $m_\Om(A):=m(A|\Om)$ .
\

It follows from \eqref{eq: ddag} that
$$e^{-sp\v}U_s^p\circ T_\Om=(U_s-1_\Om)^p=U_s^p+\sum_{k=0}^{p-1}\tbinom{p}k(-1)^kU_s^k=U_s^p+\mathcal E_{p,s}U_s^{p-1}$$
where $|\mathcal E_{p,s}|\le 2^p$. Thus
$$|U_s^p\circ T_\Om-U_s^p|\le (1-e^{-sp\v})U^p_s\circ T_\Om+2^pU_s^{p-1}$$
and
\begin{align*}
 \|U_s^p\circ T_\Om-U_s^p]|_{L^2(\Om)}& \le \|1-e^{-sp\v}\|_{L^2(\Om)}\|U^p_s\circ T_\Om\|_{L^2(\Om)}+2^p\|U_s^{p-1}\|_{L^2(\Om)}\\ &=
  \|1-e^{-sp\v}\|_{L^2(\Om)}\|U^p_s\|_{L^2(\Om)}+2^p\|U_s^{p-1}\|_{L^2(\Om)}\\ &=o(\int_\Om U_s^pdm)\ \text{as\ $s\to 0$\ by \eqref{eq: dsmathematical}}
\end{align*}
whence
$$\|\Phi_s\circ T_\Om-\Phi_s\|_{L^2(\Om)}\xrightarrow[s\to 0]{}\ 0.$$
Now suppose  that  $\Psi\in L^2(\Om),\ t_k\to 0$ so that 
\begin{align*}\Phi_{t_k}\xrightarrow[N\to\infty]{}\ \Psi\ \ \text{weakly in}\ L^2(\Om),\end{align*}
then (since $m_\Om\circ T_\Om^{-1}=m_\Om$ )
\begin{align*}\Phi_{t_k}\circ T_\Om\xrightarrow[N\to\infty]{}\ \Psi\circ T_\Om\ \ \text{weakly in}\ L^2(\Om),\end{align*}
and by  $ \|U_s^p\circ T_\Om-U_s^p]|_{L^2(\Om)}=o(\int_\Om U_s^pdm)$,
\begin{align*}\Phi_{t_k}\circ T_\Om\xrightarrow[N\to\infty]{}\ \Psi\ \ \text{weakly in}\ L^2(\Om).\end{align*}
It follows that $\Psi=\Psi\circ T_\Om$. By ergodicity, $\Psi\equiv\int_\Om\Psi dm=1$. So the only weak limit point of $\Phi_s$ as $s\to 0$ is the constant 1.
\ \ \Checkedbox{ \eqref{eq: dschemical}}

\demo{\tt Proof of \eqref{eq: dsheraldical}}

Suppose that \eqref{eq: dsheraldical} fails and let $\e>0$ and let $s_j\to 0$ be sequence so that
\begin{align*}\int_{\Om} V(p,s_j)dm\ \ {\text{\Large $\lesssim$}}\ (1-2\e)u(ps_j)\int_{\Om} V(p-1,s_j)dm.
\end{align*}
\

By \eqref{eq: a2} and Egorov's theorem, there is a subsequence  $t_k\to 0$ and $A\in\B(\Om),\ m(A)>1-\e$ so that $\widehat{U}_{t_k}\sim u(t_k)$ as $k\to\infty$ uniformly on $A$, whence
\begin{align*}\int_{\Om} V(p,t_k)dm&\ge \int_A V(p-1,t_k)\widehat{U}_{pt_k} dm
 \\ &\text{\Large $\sim$}\ u(pt_k)\int_{A} V(p-1,t_k)dm
 \\ &\overset{\eqref{eq: dschemical}}{\text{\Large $\gtrsim$}}m(A)u(pt_k)\int_{\Om} V(p-1,t_k)dm
 \\ &> (1-\e)u(pt_k)\int_{\Om} V(p-1,t_k)dm.\ \ \ \CheckedBox\text{\ \eqref{eq: dsheraldical} $\&$ \eqref{eq: dsliterary}}
\end{align*}
Next, we claim that $\exists\ M>1\ \&\ \D:\Bbb N\to\Bbb R_+$ so that
\begin{align}\label{eq: dsmilitary}\frac1M<\frac{pu(ps)}{u(s)}<M\ \ \forall\ p\ge 1,\ \ 0<s<\D(p). 
\end{align}
\pf of \eqref{eq: dsmilitary}:
\

We now use the assumption $\a=\a(T)<\infty$. Since $U_s\underset{s\to 0+}{\text{\Large $\lesssim$}}\ \a u(s)$ a.e., by Egorov's theorem,
$\exists\ A\in \B(\Om)$ with $m(A)>0$ so that $U_s\underset{s\to 0+}{\text{\Large $\lesssim$}}\ \a u(s)$ uniformly on $A$.
\

Using this, \eqref{eq: dsliterary} and \eqref{eq: dschemical}, we have
$$m(A)p!\prod_{k=1}^pu(ks)\underset{s\to 0+}{\text{\Large $\sim$}}\int_{A}U_s^pdm \underset{s\to 0+}{\text{\Large $\lesssim$}} m(A)\a^pu(s)^p.$$
Fixing $c>0$ so that 
$$(p!)^{\frac1p}\ge cp\ \forall\ p\ge 1,$$
it follows that
$$cpu(ps)\le \left(p!\prod_{k=1}^pu(ks)\right)^{\frac1p} \underset{s\to 0+}{\text{\Large $\lesssim$}}\a u(s).$$
This proves \eqref{eq: dsmilitary}.
\

Using \eqref{eq: dsmilitary}, we can now apply the de Haan-Stadtm\"uller theorem (theorem 1 in   \cite{dH-S} and theorem 2.10.2  in \cite{BGT}) that $\exists\ I>1$ so that
\begin{align*}u(s)\ \ \  =\ \ \ I^{\pm 1} a(\frac1s)
\end{align*}
thus obtaining ({\tt ER}).\ \ \ \Checkedbox

\
  \subsection*{Interarrival stochastic processes and generalized recurrent events}
 \

 Let $\xbmt$ be  a conservative, ergodic measure preserving transformation.
 \

 The {\it induced transformation} on $\Om\in\mathcal F_+$ is the probability preserving transformation
$$(\Om,\mathcal B(\Om) ,m_\Om,T_\Om)$$ where \sbul\ $m_\Om:=m(\cdot\,|\Om)$;
\sbul\ $T_\Om:\Om\to\Om$ is the {\it first return } or {\it induced} transformation  defined by \par $T_\Om x:=T^{\v(x)}x$
where $\v=\v_\Om:\Om\to\Bbb N$ is the {\it first return time} function defined by
 $\v(x):=\min\,\{n\ge 1:\ T^nx\in\Om\}$. 
 \
 
 The (one-sided)
{\it interarrival (stochastic) process} of $\Om$ is the stochastic process $(\v\circ T_\Om^n)_{n\ge 0}$ defined on $\Om$. 
It corresponds to a factor induced transformation on  $\Om$ corresponding to
the sub-invariant factor algebra $\B_0:=\s(\{T^{-n}\Om:\ n\ge 0\})$.

\

As in \cite{AN}, a stochastic  process $(X_1,X_2,\dots)$ is {\it continued fraction mixing}  if
$\vartheta(1)<\infty\ \&\  \vartheta(n)\downarrow 0$ where
$$\vartheta(n):=\sup\{|\tfrac{\Bbb P(A\cap B)}{\Bbb P(A)\Bbb P(B)}-1|:\
A\in\s_{1}^k,\ B\in\s_{k+n}^\infty,\ \Bbb P(A)\Bbb P(B)>0,\ k\ge
1\}.$$ Here, $\s_k^N$ denotes the $\s$-algebra generated by the random variables
$\{X_j: k\le j<N+1\}$ for $k<N+1\le\infty$.

\

Let $\xbmt$ be a conservative, ergodic measure preserving transformation.
\

Let $\Om\in\mathcal F_+$. Consider the property: 
\f (\Bicycle)  \ \ $\Om$'s interarrival stochastic process $(\v_\Om\circ T^n)_{n\ge 0}$ is 
continued fraction  mixing with coefficients satisfying $\sum_{n=1}^\infty\frac{\vartheta(n)}n<\infty$.
\

We'll call any $\Om\in\mathcal F_+$ satisfying (\Bicycle) a  (\Bicycle) {\it set }. 
\subsection*{Remarks}
\

\sms (i) Any recurent event (as in 5.2 of \cite{IET}) has an independent, interarrival stochastic process whence is a  (\Bicycle) set. 

\sms (ii)\ Examples are also obtained by noting that (as shown in \cite{GM}) any stationary stochastic process driven by a mixing 
Gibbs-Markov map and with 
observable measurable with respect to the Markov partition is continued fraction  mixing with exponentially decaying coefficients.

\sms (iii) By lemma 3.7.4 in  \cite{IET}, a transformation with a  (\Bicycle) set   has a factor 
where the  (\Bicycle) set is a Darling-Kac set (and is hence  pointwise dual ergodic). 

\

\

For   $\Om\in\mathcal F_+$,  set
 $$L(t)=L_\Om(t):=\int_\Om(\v_\Om\wedge t)dm_\Om,\ \ \ a(t):=\frac{t}{L(t)}\ \ \&\ \ b:=a^{-1}.$$

\proclaim{Theorem 5}
\

Let $\xbmt$ be a conservative, ergodic measure preserving transformation equipped with a  (\Bicycle) set $\Om\in\mathcal F_+$, then $\xbmt$ satisfies {\Large\Pointinghand}   if and only if
 \begin{align}\label{eq: symqueen}\sum_{n=1}^\infty n\cdot\left(\frac{m_\Om([\v_\Om\ge n])}{L_\Om(n)}\right)^2<\infty.
 \end{align}

In this case
\begin{align}\label{eq: symking}\b(T)=1\ \&\ \u\b(T)=\frac12.
\end{align}\endproclaim
\subsection*{Remarks} 
\

1.\ \ As shown in \cite{AN} 
the condition\  \eqref{eq: symqueen} \ characterizes the ``trimmed sum" convergence property:
$$\frac1{b(n)}(\v_n-\max_{0\le k\le n-1}\v\circ T_\Om^k)\ \xyr[n\to\infty]{}\ 1\ \text{a.s.}$$
where $\varphi_n:=\sum_{k=0}^{n}\varphi_\Omega\circ T_\Omega^k$. See the earlier \cite{M1, M2} for the independent case and \cite{DV} for the case of continued fraction partial quotients. 
\

2.  By lemma 3 in \cite{M2}, the condition \eqref{eq: symqueen} implies that $L_\Om(n)$ is slowly varying.
 
\demo{Proof}
\

 Define  $\v^{\pm}:\Om\to\Bbb N$ by
$$\v^{\pm}(x):=\inf\,\{n\ge 1:\ T_\Om^{\pm n}(x)\in\Om\}\ \&\ \v^{\pm}_J:=\sum_{j=0}^{J-1}\v^{\pm}\circ T_\Om^{\pm j}$$ and define
$$\v_n(x):=\begin{cases} &\v^{+}_n(x)\ \ \ \ \ \ \ \ \ n\ge 1,\\ & 0\ \ \ \ \ \ \ \ \ n=0,\\ & -\v_{-n}^{-}(\s^{-1}x)\ \ \ \ \ \ \ \ \ n\le -1.\end{cases}$$
It follows that
\begin{align}\label{eq: dsrailways}\Si_n(1_\Om)(x)=\#\{k\in\Bbb Z:\ |\v_k(x)|\le n\}.
 \end{align}

For $n\in\Bbb N,\ t>0$:
Define $B_n(t)$ by
$$B_n(t):=\bigcup_{k=2^n+1}^{2^{n+1}}[\v^{-}\circ T_\Om^{-k}>tb(2^n)]\cap [\v^{+}\circ T_\Om^{k}>tb(2^n)] $$  where $b=a^{-1},\ a(n)=a_n(T).$

\

Following the ideas in the proof of lemma 1.2 in \cite{AN}, we claim that
\begin{align}\label{eq: Leftscissors}P(B_n(t))\asymp 2^nP([\v^{+}>tb(2^n)])^2
\end{align}
\demo{Proof of \eqref{eq: Leftscissors}}\ \
Evidently,
\begin{align*}P(B_n(t))&\le\sum_{k=2^n+1}^{2^{n+1}}P([\v^{-}\circ T_\Om^{-k}>tb(2^n)]\cap [\v^{+}\circ T_\Om^{k}>tb(2^n)])\\ &\le\sum_{k=2^n+1}^{2^{n+1}}(1+\vartheta(2^{n+1}))P([\v^{-}\circ T_\Om^{-k}>tb(2^n)])P( [\v^{+}\circ T_\Om^{k}>tb(2^n)])\\ &=
(1+\vartheta(2^{n+1}))2^nP([\v^{+}>tb(2^n)])^2.
\end{align*}
For the other inequality, choose $\kappa\ge 1$ so that $\vartheta(\kappa)<\tfrac12\ \&$ \f $2^{n+1}P([\v^{\pm}>tb(2^n)])\le\tfrac12\ \forall\ n\ge\kappa$.
\

 Fix $n\gg \kappa\ \&\ 2^n< k\le 2^{n+1}$ and define
$$A_k^{(n)}:=[\v^{+}\circ T_\Om^k\wedge\v^{-}\circ T_\Om^{-k}>tb(2^n)]\cap \bigcap_{2^n< j\le 2^{n+1},\ |j-k|\ge\kappa}[\v^{+}\circ T_\Om^j\wedge\v^{-}\circ T_\Om^{-j}\le tb(2^n)].$$
It follows that
\begin{align*}P(&A_k^{(n)})=\\ &P([\v^{+}\circ T_\Om^k\wedge\v^{-}\circ T_\Om^{-k}>tb(2^n)]\cap \bigcap_{2^n< j\le 2^{n+1},\ |j-k|\ge\kappa}[\v^{+}\circ T_\Om^j\wedge\v^{-}\circ T_\Om^{-j}\le tb(2^n)])\\ &\ge (1-\vartheta(\kappa))^3P([\v>tb(2^n])^2P(\bigcap_{2^n< j\le 2^{n+1}}[\v^{+}\circ T_\Om^j\wedge\v^{-}\circ T_\Om^{-j}\le tb(2^n)])^3\\ &=\tfrac18P([\v>tb(2^n])^2(1-P(\bigcap_{2^n< j\le 2^{n+1}}[\v^{+}\circ T_\Om^j\vee\v^{-}\circ T_\Om^{-j}\le tb(2^n)]))^3\\ &\ge\tfrac18P([\v>tb(2^n])^2(1-2^{n+1}P([\v>tb(2^n)]))^3\\ &\ge
\frac1{64}P([\v>tb(2^n])^2.
\end{align*}
Moreover
$$\sum_{2^n< k\le 2^{n+1}}1_{A_k^{(n)}}\le\ (2\kappa+1) 1_{B_n(t)}$$
whence
\begin{align*}P(B_n(t)) &\ge\frac1{2\kappa+1}\sum_{2^n< k\le 2^{n+1}}P(A_k^{(n)})\\ &\ge
\frac{2^n}{64(2\kappa+1)}P([\v>tb(2^n])^2.\ \ \CheckedBox\ \text{\eqref{eq: Leftscissors}}
\end{align*}
It follows from \eqref{eq: Leftscissors} and continued fraction mixing that
$$P(B_n(t)\cap B_{n'}(t))\ \asymp\ P(B_n(t))P(B_{n'}(t))\ \ \text{for}\ n\ne n'\in \Bbb N.$$

\

The Borel Cantelli lemmas now ensure (as in \cite{AN} $\&$ \cite{M1, M2}) that
\begin{equation} \label{eq: dstechnical}
\begin{aligned}
\sum_{n= 1}^\infty 1_{B_n(t)}=\infty\ \text{a.s.}& \iff\ \sum_{n= 1}^\infty 2^nP([\v^{+}>tb(2^n)])^2  =\infty\\ & \ \iff\ \sum_{n= 1}^\infty P([\v^{+}>tb(n)])^2 =\infty.
\end{aligned}
\end{equation}

If, in addition, $b$ is {\it weakly regularly varying} in the sense that
$$\exists\ M>1\ \text{such that}\ A(2t)\le MA(t)\ \&\ 2A(t)\le A(Mt)\ \forall\ \text{large}\ t\in\Bbb R_+,$$
then
the convergence of \eqref{eq: dstechnical}\  for some $t>0$ implies its convergence  for every $t>0$; a situation  characterized by  \eqref{eq: symqueen}\  (for more details, see \cite{AN}).
\

\

To continue, we pass to  the {\it  one-sided factor} $$\pi:\xbmt\to (X_0,\B_0,m_0,T_0)$$ defined by
$$\pi^{-1}\B_0=\mathcal F_\Om=\mathcal F_\Om:=\s(\{T^{-n}\Om:\ n\ge 0\}).$$
 Fix $\Om_0\in\B_0,\ \pi^{-1}\Om_0=\Om$, then
 $\Om_0$ is a Darling Kac set for $T_0$.
 
\

 \demo{Proof that \eqref{eq: symqueen}\  $\Rightarrow$  {\Large\Pointinghand} $\&$ \eqref{eq: symking}}
 \

\

Suppose that \eqref{eq: symqueen} \ is satisfied then, as above, 
$L(n)$ is slowly varying and  by the asymptotic renewal equation (3.8.6 in \cite{IET}) 
$a_n(T)\propto\frac{n}{L(n)}$ is $1$-regularly varying and, in particular, weakly regularly varying.
\

Moreover,  by \cite{AN},  $\varliminf_{n\to\infty}\frac{\v_n}{b(n)}=1$ a.s., whence (see \cite{LB}) $\a(T)=1$.
\

By proposition 1,
$$\varliminf_{n\to\infty}\frac1{2a_n(T)}{\Si_n(1_A)}\le \frac{m(A)}2\ \text{a.e.}\  \forall\ A\in\B,\ 0<m(A)<\infty.$$
Next, $\sum_{k=1}^\infty 1_{B_k}<\infty$ a.s. and by theorem 1.1 in \cite{AN}, $\exists\ \e:\Bbb N\x\Om\to\{-,+\}$ so that
$$\frac{\v^{\e(n,x)}_n(x)}{b(n)}\xyr[n\to\infty]{}\ 1\ \text{a.s.  where}\ b=a^{-1},\ a(n)=a_n(T).$$
In addition, it follows that a.s.,
\begin{align*}&\frac1{a_n(T)}S_n^{(T^{\e(a(n),x)})}(1_\Om)(x)\ \ \xyr[n\to\infty]{}\ 1\ \text{a.s., whence}\\ &   \varliminf_{n\to\infty}\frac1{2a_n(T)}{\Si_n(1_\Om)}\ge \frac12.\ \ \ \CheckedBox\text{{\Large\Pointinghand} $\&$ \eqref{eq: symking}}
\end{align*}

\demo{Proof that\ \ {\Large\Pointinghand} $\Rightarrow$  \eqref{eq: symqueen}}
\

It follows from $\a(T)<\infty$ via Theorem 4 that  $b$ is weakly regularly varying.

If \eqref{eq: symqueen}\  fails, then as above, for every $t>0$
$\sum_{k=1}^\infty 1_{B_k(t)}=\infty$ a.s. and a.e. $x\in\Om,\ \exists\ n_k\to\infty$ so that
$$\v_{n_k}^{+}(x),\ \v_{n_k}^{-}(x)\ >tb(n_k).$$

Set $N_k:=\lfl tb(n_k)\rfl$, then
$$S_{N_k}^{(T^{\pm 1})}(1_\Om)<n_k\sim a(\frac{N_k}t)\le \frac{MI^2}t a(N_k).$$
It follows that
\begin{align*}
 \varliminf_{n\to\infty}\frac1{2a(N)}{\Si_n(1_\Om)}(x)&\le
 \varlimsup_{k\to\infty}\frac1{2a(N_k)}{\Si_{N_k}(1_\Om)}(x)\\ &\le
 \varlimsup_{k\to\infty}\frac1{2a(N_k)}{S_{N_k}^{(T)}(1_\Om)}(x)+
 \varlimsup_{k\to\infty}\frac1{2a(N_k)}{S_{N_k}^{(T^{- 1})}(1_\Om)}(x)\\ &\le
 \frac{MI^2}t\xyr[t\to\infty]{}\ 0.\ \ \CheckedBox
\end{align*}

\section*{\S5 The multidimensional situation}
\

\subsection*{Example 5.1}
\

Let $\xbm$ be $\Bbb R$ equipped with Borel sets and Lebesgue measure.
Let $\a,\ \b\in\Bbb R$ be linearly independent over $\Bbb Q$ and define
$$\tau =\tau^{(\a,\b)}:\Bbb Z^2\to\text{\tt MPT}\,\xbm$$
by
$$\tau_{(k,\ell)}(x):=x+k\a+\ell\b.$$
Define
$$\Xi_n^{(\tau)}(f):=\sum_{|k|,\ |\ell|\le N}f\circ \tau_{(k,\ell)}.$$
We claim that
\begin{align*}\tag{{\leftthumbsup}}\frac{\Xi_n^{(T)}(f)}{2N+1}\
\xyr[N\to\infty]{}\ R\int_Xfdm\ \ \text{a.e.}\ \forall\ f\in\ L^1(m)
\end{align*}
where $R:=\frac{\min\,\{|\a|,|\b|\}}{\max\,\{|\a|,|\b|\}}$.
\demo{Proof of ({\leftthumbsup}) when $|\a|>|\b|=1$}
\

Here $R=\frac1{|\a|}\ \ \&\a\notin\Bbb Q$. We have that $W=[0,1)$ is a {\it maximal wandering set} for $\tau_{0,1}$ in the sense that
$$X=\bigcupdot_{n\in\Bbb Z}\tau_{0,n}W,$$
whence, since $|\a|>1$,  $\exists\ \kappa:\Bbb Z\x W\to\Bbb Z$ so that for $x\in W$,
$$\frak n(x,W):=\{u\in\Bbb Z^2:\ \tau_u(x)\in W\}=\{(\ell,\kappa(\ell,x)):\ \ell\in\Bbb Z\}.$$
Here $|\kappa(\ell,x)|=|\ell\a|\pm 1\ \forall\ \ell\in\Bbb Z$, whence for $N\ge 1,\ x\in W$,
$$\frak n(x,W)\cap[-N,N]^2=\{(\ell,\kappa(\ell,x)),\ell):\ \ell\in\Bbb Z,\ |\ell|,\ |\kappa(\ell,x)|\le N\}$$
and
$$\Xi_N^{(\tau)}(1_W)(x)=\#\{(\kappa(\ell,x),\ell):\ \ell\in\Bbb Z,\ |\ell|\le N\}\sim\frac{ 2N}{|\a |}.$$
Next define $S:W\to W$ by $S(x):=\tau_{(1,\kappa(1,x))}$, then
$$S(x)=x+\a\ \mod\ 1.$$
Thus $\tau$ is ergodic and for $f:X\to\Bbb R$, supported and continuous on $W$, we have on $W$:
\begin{align*}\frac{\Xi_n^{(T)}(f)}{2N+1}&=\frac1{2N+1}\sum_{|k|,\ |\ell|\le N}f\circ \tau_{(k,\ell)}
\\ &\sim\frac1{2N}\sum_{|\ell|\le |\a|N}f\circ S^\ell\\ &\xyr[N\to\infty]{}\int_Wfdm_W=\frac1{|\a|}\int_Xfdm\ \ \text{uniformly on}\ W.
\end{align*}

The proposition follows from this via \cite{H}.\ \ \Checkedbox

\

It is not hard to show that
\sbul the above action $T$ is uniquely ergodic in the sense that the only $T$-invariant Radon measures on $\Bbb R$ are multiples of $m$; and
\sbul the convergence ({\leftthumbsup}) is uniform on compact subsets for bounded continuous functions $f$.
\subsection*{Example 5.2}
\

Let $f\in\mathcal P(\Bbb N)$, and let $\Om:=\Bbb N^\Bbb Z$ and $P=P_f\in\mathcal P(\Om)$ be product measure defined by
$$P(\{\om\in\Om:\ \om_{k+i}=n_i\ \forall\ 1\le i\le N\})=\prod_{1\le i\le N}f_{n_i}\ \ \ \ \ (k\in\Bbb Z).$$
Let
$$\xbm:=(\Om\x\Bbb Z,\B(\Om\x\Bbb Z),P_f\x\#),$$ let $\s:\Om\to\Om$ is the shift and $$\psi:\Bbb Z^2\to\text{\tt MPT}\,\xbm$$
by
$$\psi_{1,0}(\om,n):=(\s\om,n+\om_0)\ \ \&\  \psi_{0,1}(\om,n):=(\om,n+1).$$
The action is ergodic since for $\om\in\Om$,
$$\{\psi_{k,\ell}(\om,0):\ k,\ \ell\in\Bbb Z\}=\bigcup_{n\in\Bbb Z}\{\s^j(\om):\ j\in\Bbb Z\}\x\{n\}.$$
Moreover, writing
$$s_k(\om):=\begin{cases} &\sum_{j=0}^{k-1}\om_j\ \ \ \ \ k\ge 1,\\ &
0\ \ \ \ \ \ \ \ \ \ k=0,\\ & -\sum_{j=1}^{|k|}\om_{-j}=-s_{-k}(\s^k\om)\ \ \ \ \ k\le -1,;\end{cases}$$
we have
\begin{align*}\sum_{|k|,\ |\ell|\le N}1_{\Om\x\{0\}}(\psi_{k,\ell}(\om,0))&=\sum_{|k|,\ |\ell|\le N}1_{\Om\x\{0\}}((\s^k\om,s_k(\om)-\ell))\\ &=
\#\{k\in [-N,N]:\ |s_k(\om)|\le N\}.
\end{align*}
Let $u=(u_0,u_1,\dots)$ be the renewal sequence with lifetime distribution $f$ and let $a_u(n):=\sum_{k=1}^nu_k$.

\

 By \eqref{eq: dsrailways} and theorem 5, we have that  the following conditions (on $f\in\mathcal P(\Bbb N)$) are equivalent:
\begin{align*}&\tag{B1}\exists\ a_n>0\ \ \text{so that}\ \ \Xi^{(\psi)}(1_{\Om\x\{0\}})\ \asymp\ a_n
\\ &\tag{B2}
\Xi^{(\psi)}(1_{\Om\x\{0\}})\ \asymp\ a_u(n);\\ &\tag{B3}\sum_{n=1}^\infty\left(\frac{f([n,\infty))}{L_f(n)}\right)^2<\infty\ \text{where}\ L_f(n):=\sum_{k=1}^nf([k,\infty)).
\end{align*}

In this case (\leftthumbsup) fails.

\

 The above examples show that a conservative, ergodic, infinite measure preserving $\Bbb Z^2$ action  having a dissipative generator with a maximal wandering set of finite measure can
 \sbul satisfy (\leftthumbsup);
 \sbul satisfy   (\Pointinghand) while not satisfying (\leftthumbsup),
 \sbul not satisfy  (\Pointinghand).

\

It follows from theorem 2 that a a conservative, ergodic, infinite measure preserving $\Bbb Z^2$ action  having a dissipative generator with a maximal wandering set of infinite measure cannot
  satisfy (\leftthumbsup), the other two possibilities being available.

\

\subsection*{Question}\ \ \ There are conservative, ergodic, infinite measure preserving $\Bbb Z^2$ actions with both generators conservative.
We do not know which of the above possibilities are available for such an action.


\begin{bibdiv}
\begin{biblist}

\bib{BRE}{article}{
      author={Aaronson, J.},
       title={Rational ergodicity, bounded rational ergodicity and some
  continuous measures on the circle},
        date={1979},
        ISSN={0021-2172},
     journal={Israel J. Math.},
      volume={33},
      number={3-4},
       pages={181\ndash 197 (1980)},
         url={http://dx.doi.org/10.1007/BF02762160},
        note={A collection of invited papers on ergodic theory},
      review={\MR{571529 (81f:28012)}},
}

\bib{IET}{book}{
      author={Aaronson, Jon},
       title={An introduction to infinite ergodic theory},
      series={Mathematical Surveys and Monographs},
   publisher={American Mathematical Society},
     address={Providence, RI},
        date={1997},
      volume={50},
        ISBN={0-8218-0494-4},
      review={\MR{1450400 (99d:28025)}},
}

\bib{LB}{incollection}{
      author={Aaronson, Jon},
      author={Denker, Manfred},
       title={Lower bounds for partial sums of certain positive stationary
  processes},
        date={1989},
   booktitle={Almost everywhere convergence ({C}olumbus, {OH}, 1988)},
   publisher={Academic Press},
     address={Boston, MA},
       pages={1\ndash 9},
      review={\MR{1035234 (91e:28012)}},
}

\bib{GM}{article}{
      author={Aaronson, Jon.},
      author={Denker, Manfred},
       title={Local limit theorems for partial sums of stationary sequences
  generated by {G}ibbs-{M}arkov maps},
        date={2001},
        ISSN={0219-4937},
     journal={Stoch. Dyn.},
      volume={1},
      number={2},
       pages={193\ndash 237},
         url={http://dx.doi.org/10.1142/S0219493701000114},
      review={\MR{1840194 (2002h:37014)}},
}

\bib{AN}{article}{
      author={Aaronson, Jon},
      author={Nakada, Hitoshi},
       title={Trimmed sums for non-negative, mixing stationary processes},
        date={2003},
        ISSN={0304-4149},
     journal={Stochastic Process. Appl.},
      volume={104},
      number={2},
       pages={173\ndash 192},
         url={http://dx.doi.org/10.1016/S0304-4149(02)00236-3},
      review={\MR{1961618 (2004c:60077)}},
}

\bib{AZ}{article}{
      author={Aaronson, Jon},
      author={Zweimueller, Roland},
       title={Limit theory for some positive, stationary processes with
  infinite mean},
        date={2010},
      eprint={arXiv/1008.3919},
}

\bib{BGT}{book}{
      author={Bingham, N.~H.},
      author={Goldie, C.~M.},
      author={Teugels, J.~L.},
       title={Regular variation},
      series={Encyclopedia of Mathematics and its Applications},
   publisher={Cambridge University Press},
     address={Cambridge},
        date={1987},
      volume={27},
        ISBN={0-521-30787-2},
      review={\MR{898871 (88i:26004)}},
}

\bib{CH}{incollection}{
      author={Chacon, R.~V.},
       title={A geometric construction of measure preserving transformations},
        date={1967},
   booktitle={Proc. {F}ifth {B}erkeley {S}ympos. {M}ath. {S}tatist. and
  {P}robability ({B}erkeley, {C}alif., 1965/66), {V}ol. {II}: {C}ontributions
  to {P}robability {T}heory, {P}art 2},
   publisher={Univ. California Press},
     address={Berkeley, Calif.},
       pages={335\ndash 360},
      review={\MR{0212158 (35 \#3033)}},
}

\bib{DGPS}{article}{
      author={Dai, Irving},
      author={Garcia, Xavier},
      author={Puadurariu, Tudor},
      author={Silva, Cesar~E.},
       title={On rationally ergodic and rationally weakly mixing rank-one
  transformations},
        date={2012},
      eprint={arXiv/1208.3161},
}

\bib{dH-S}{article}{
      author={de~Haan, L.},
      author={Stadtm{\"u}ller, U.},
       title={Dominated variation and related concepts and {T}auberian theorems
  for {L}aplace transforms},
        date={1985},
        ISSN={0022-247X},
     journal={J. Math. Anal. Appl.},
      volume={108},
      number={2},
       pages={344\ndash 365},
         url={http://dx.doi.org/10.1016/0022-247X(85)90030-7},
      review={\MR{793651 (86k:44002)}},
}

\bib{DV}{article}{
      author={Diamond, Harold~G.},
      author={Vaaler, Jeffrey~D.},
       title={Estimates for partial sums of continued fraction partial
  quotients},
        date={1986},
        ISSN={0030-8730},
     journal={Pacific J. Math.},
      volume={122},
      number={1},
       pages={73\ndash 82},
         url={http://projecteuclid.org/euclid.pjm/1102702122},
      review={\MR{825224 (87f:11056)}},
}

\bib{DE}{article}{
      author={Dowker, Yael~Naim},
      author={Erd{\H{o}}s, Paul},
       title={Some examples in ergodic theory},
        date={1959},
        ISSN={0024-6115},
     journal={Proc. London Math. Soc. (3)},
      volume={9},
       pages={227\ndash 241},
      review={\MR{0102584 (21 \#1374)}},
}

\bib{FR-Book}{book}{
      author={Friedman, Nathaniel~A.},
       title={Introduction to ergodic theory},
   publisher={Van Nostrand Reinhold Co.},
     address={New York},
        date={1970},
        note={Van Nostrand Reinhold Mathematical Studies, No. 29},
      review={\MR{0435350 (55 \#8310)}},
}

\bib{FR-AMM}{article}{
      author={Friedman, Nathaniel~A.},
       title={Replication and stacking in ergodic theory},
        date={1992},
        ISSN={0002-9890},
     journal={Amer. Math. Monthly},
      volume={99},
      number={1},
       pages={31\ndash 41},
         url={http://dx.doi.org/10.2307/2324545},
      review={\MR{1140275 (93d:28022)}},
}

\bib{GN}{article}{
      author={Gorodnik, Alex},
      author={Nevo, Amos},
       title={Ergodic theory and the duality principle on homogeneous spaces},
        date={2012},
      eprint={arXiv/1205.4413},
}

\bib{H}{article}{
      author={Hochman, Michael},
       title={A ratio ergodic theorem for multiparameter non-singular actions},
        date={2010},
        ISSN={1435-9855},
     journal={J. Eur. Math. Soc. (JEMS)},
      volume={12},
      number={2},
       pages={365\ndash 383},
         url={http://dx.doi.org/10.4171/JEMS/201},
      review={\MR{2608944 (2011g:37006)}},
}

\bib{EH}{book}{
      author={Hopf, E.},
       title={Ergodentheorie},
      series={Ergebnisse der Mathematik und ihrer Grenzgebiete, 5. Bd},
   publisher={Julius Springer},
        date={1937},
      number={v. 5, no. 2},
         url={http://books.google.co.il/books?id=icHQpwAACAAJ},
}

\bib{MS}{techreport}{
      author={Maucourant, François},
      author={Schapira, Barbara},
       title={Distribution of orbits in r$^2$ of a finitely generated group of
  sl(2,r),},
        date={2012},
      number={arXiv:1204.5158},
}

\bib{M1}{article}{
      author={Mori, Toshio},
       title={The strong law of large numbers when extreme terms are excluded
  from sums},
        date={1976},
     journal={Z. Wahrscheinlichkeitstheorie und Verw. Gebiete},
      volume={36},
      number={3},
       pages={189\ndash 194},
      review={\MR{0423494 (54 \#11470)}},
}

\bib{M2}{article}{
      author={Mori, Toshio},
       title={Stability for sums of i.i.d. random variables when extreme terms
  are excluded},
        date={1977},
     journal={Z. Wahrscheinlichkeitstheorie und Verw. Gebiete},
      volume={40},
      number={2},
       pages={159\ndash 167},
      review={\MR{0458542 (56 \#16742)}},
}

\bib{N1}{book}{
      author={Nadkarni, M.~G.},
       title={Spectral theory of dynamical systems},
      series={Texts and Readings in Mathematics},
   publisher={Hindustan Book Agency},
     address={New Delhi},
        date={2011},
      volume={15},
        ISBN={978-93-80250-21-2},
        note={Reprint of the 1998 original},
      review={\MR{2847984 (2012h:37005)}},
}

\bib{N2}{book}{
      author={Nadkarni, M.~G.},
       title={Basic ergodic theory},
     edition={Third},
      series={Texts and Readings in Mathematics},
   publisher={Hindustan Book Agency},
     address={New Delhi},
        date={2013},
      volume={6},
        ISBN={978-93-80250-43-4},
      review={\MR{2963410}},
}

\bib{CESAR}{book}{
      author={Silva, C.~E.},
       title={Invitation to ergodic theory},
      series={Student Mathematical Library},
   publisher={American Mathematical Society},
     address={Providence, RI},
        date={2008},
      volume={42},
        ISBN={978-0-8218-4420-5},
      review={\MR{2371216 (2009d:37001)}},
}

\end{biblist}
\end{bibdiv}
\end{document}